\documentstyle[11pt]{article}
\title{The variance of   the incipient   infinite cluster in two-dimensional  percolation
 }
\author{Yu Zhang}
\begin{document}
\baselineskip .20in
\maketitle
\begin{abstract}
Consider bond percolation on the square lattice.  Let ${\cal C}$ be the incipient infinite cluster with the incipient  measure $\nu$.
If a one-arm path  exponent exists and equals  $5/48$,
it is well known that
 $E_\nu|{\cal C}\cup [-n,n]^2| = n^{ 91/48+o(1)}$.
In this paper, we focus on the variance of $|{\cal C}\cup [-n,n]^2|$ and show that
$\sigma^2_\nu (|{\cal C}\cup [-n,n]^2|)= n^{91/24+o(1)}$.

\end{abstract}

\section{ Introduction and statement of results.}
Consider bond percolation on the square lattice ${\bf Z}^2$.
 For a given $(x, y)\in {\bf Z}^2$,
its nearest four neighbors are defined as $(x\pm 1, y)$ and $(x, y\pm1)$.
Edges between neighboring or adjacent vertices therefore correspond to vertical or horizontal displacements of one unit.
Each edge  is independently {\em open} with 
probability $p$ or {\em closed} with probability $1-p$. 
The corresponding probability measure on the configurations of open
and closed edges is denoted by $P_p$. We also denote by $E_p$  and $\sigma_p$ the expectation and the variance with respect to $P_p$. 
A path from $u$ to $v$ is a sequence $(v_0,e_1, v_1,... ,v_{i-1}, e_i,  v_{i},... ,v_n)$
with adjacent  vertices $v_i$  and $v_{i+1}$  by edge $e_i$ such that $v_0=u$ and $v_n=v$. 
A circuit is a path with distinct vertices $v_i$ ($1\leq i\leq n-1$) except $v_0=v_n$. 
A path is called open or closed if all of its edges  are open or closed. The occupied cluster of the vertex $x$, 
${\bf C}(x)$,
consists of all open edges that are connected to $x$ by an open path.
For any collection $A$ of vertices, $|A|$ 
denotes the cardinality of $A$. We choose ${\bf 0}$ as the origin.

The {\em percolation probability} is
\begin{eqnarray*}
\theta (p)= P_p(|{\bf C}({\bf 0})|=\infty),
\end{eqnarray*}
and the {\em critical probability} is 
$$p_c=\sup\{p:\theta (p)=0\}.$$
It is well known   (see chapter 3 in Kesten (1982)) that for bond percolation on the square lattice 
$$p_c=0.5.\eqno{(1.1)}$$

It is also well known that there is no infinite open cluster at $p_c$. However, the mean open cluster size is infinity at $p_c$. 
In other words, there should be a large open cluster around the origin at $p_c$. Kesten (1986) proposed  building an infinite open cluster by conditioning on the existence of a large cluster at the origin. More precisely, for any cylinder event
${\cal A}$, we write
$$\nu_n({\cal A})=P_{p_c} ({\cal A}\,\,\,|\,\,\, {\bf 0} \rightarrow \partial [-n,n]^2),$$
where $A\rightarrow B$ means that there is an open path from $A$ to $B$.  
Kesten (1986) proved that
$$ \lim_{n\rightarrow \infty} \nu_n(A) = \nu (A) \mbox{ for a cylinder event }A.\eqno{(1.2)}$$
The probability measure $\nu$ is called the {\em incipient  measure}.
Let $E_\nu(\cdot) $ and $\sigma _\nu^2(\cdot)$ be the expected value and the variance respect with measure to  $\nu$.
It is clear that the origin has $\nu$-probability one of being in an infinite open cluster  ${\cal C}$.
Physicists are interested in the geometry of the infinite incipient  cluster.  They believe that ${\cal C}$ should have a fractal between one and two. 
Before introducing the geometry results for
the infinite incipient cluster, we would like to introduce a few basic estimates in percolation.
Let us introduce $k$-arm paths. Let $B(n)=[-n,n]^2$ and $B_e(n)$ be  all the edges in $[-n,n]^2$. Consider the annulus 
$$A(m,n)=\{B(n)\setminus B(m)\}\cup \{\partial B(m)\}\mbox{ for } m <n.$$
Let ${\cal Q}_k(b, m,n)$ be the event that there exist $i$ disjoint occupied paths and $j$ disjoint vacant paths with $i+j=k$ for all $i,j\geq 1$
from $b+\partial B(m)$ to $b+\partial B(n)$ inside $A(m,n)$.
We call them $k$-{\em arm paths}.
For simplicity, let ${\cal Q}_k(m,n)={\cal Q}_k({\bf 0}, m, n)$ and  ${\cal Q}_k(n)={\cal Q}_k({\bf 0},0, n)$.
It is believed by Aizenman, Dulpantier, and Alharony (1999) that
$$P_{p_c}({\cal Q}_k(m,n))= \left({m\over n}\right)^{(k^2-1)/12+o(1)}\eqno{(1.3)}$$
 for $k \geq 2$ and
$$P_{p_c}({\cal Q}_1(m,n))= \left({m\over n}\right)^{5/48+o(1)}\eqno{(1.4)}$$
for a fixed $m$ as $n\rightarrow \infty$.
In fact, it is more important to show (1.3) when  $k=4$. Fortunately, 
by using the Schramm-Loewner evolution (SLE) argument and 
 Smirnov's  scaling limit on the triangular lattice,
(1.3) and (1.4)  were  proved  (see Lawler, Schramm, and Werner  (2002) for $k=1$ and  Smirnov and Werner's Theorem 4  (2001) for  $k=4$). 
Let $T_n=|{\cal C}\cap [-n,n]^2|$. One of most important questions is to ask what the fractal of $T_n$ is.
 With (1.4), Kesten (1987) showed that 
$$ E_\nu T_n=n^{2-5/48+o(1)}= n^{91/48+o(1)}.\eqno{(1.5)}$$
In this paper, we investigate the variance of $T_n$ to show the following result.\\

{\bf Theorem.}  {\em For each $n$, if (1.4) holds, then
$$\sigma^2_\nu (T_n )= n^{91/24+o(1)}.$$}

{\bf Remarks.} 1.  By Kesten's argument in (1986) (see (2.1) below), we know that
$$ E_\nu T_n^2 = n^{91/24+o(1)}.\eqno{(1.6)}$$
Thus, $\sigma_\nu^2( T_n )\approx  E_\nu T_n^2$.  This tells us that  $T_n$ has a very large variance as the same order as its second moments. In other words, $T_n$ has a very large tail.\\
2.  We may define $S_n$ be the number of vertices in $B(n)$ connecting to 
$\partial B(2n)$ by open paths.  By (1.4), it is easy to obtain that
$$ E_{p_c}(S_n) = n^{2-5/48+o(1)}.$$
Similarly, by  Kesten's argument above, it can be show that
$$E_{p_c}S_n^2=n^{4-5/24+o(1)}=n^{91/24+o(1)}.$$
We want to mention that by using the same proof of Theorem 1, we can also show that
$$\sigma_{p_c}^2( S_n) =n^{4-5/24+o(1)}=n^{91/24+o(1)}.\eqno{(1.7)}$$
If $ p> p_c$, then there is an infinite open cluster  from the origin with a positive probability.  On the existence of this cluster, we may still denote by $T_n$ the open cluster from the origin inside $B(n)$. By Theorem 4 in Zhang (2001),
$$\sigma_p^2(T_n)\approx n.$$
Furthermore,  $T_n$ satisfies the central limit theorem (CLT) if $p > p_c$. However,  since $T_n$ has a large tail, we  do not think  that $T_n$ satisfies the CLT at $p_c$.\\
3. For the triangular lattice, since (1.4) holds as we mentioned above, the result of the theorem will hold for the triangular lattice. Similarly, if (1.4) holds for the other two-dimensional lattices, we can also show the theorem for these lattices.

\section{ Preliminaries.}  

In this section, we  introduce a few basic properties  and estimates of bond percolation in the square lattice.  Most results are obtained from
Kesten (1982) and (1987).
For any $u,v\in {\bf Z}^2$, let $d(u,v)$ be the Euclidean distance between $u$ and $v$.
 We define the dual lattice ${{\bf Z}^*}^2=\{v+(0.5,0.5):v\in {\bf Z}^2\}$ and edges joining all pairs of vertices, which are a unit distance apart.
 For each edge set $A$, we may denote by $A^*$ its dual edges.
 For example, $e^*$ is bisected by the edge $e$.
 Given a finite connected graph $G$,   a vertex $u\not\in G$, but is adjacent
to $G$, is called the {\em boundary} of $G$. We denote by $\partial G$  the boundary vertices of $G$.
If a vertex $v\in \Delta G$, and there is an infinite path from $v$ without using $G$, then $v$ is called  the {\em exterior
boundary}. We denote by $\Delta G$  the exterior boundary of $G$. 
The edges not in $G$, but adjacent to $\Delta G$, are called the exterior  edges.
We introduce
a topology result (see Lemma 2.23 in Kesten (1982) or Proposition 11.2 in Grimmett (1999)).\\

{\bf Lemma 2.1.} {\em If $G$ is a finite cluster, then 
$\Delta  G$ is a circuit containing $G$ in its interior.  Furthermore,  $\Delta G$ is the smallest 
vacant circuit containing $G$ in its interior if $G$ is occupied.}\\

At $p_c$, by  symmetry, we can show that there is an open crossing or a closed dual crossing  in any  square with a positive probability at $p_c$. We can also extend this crossing
two times  longer  with a positive probability by using the RSW lemma. Together with the FKG inequality,  we can show that there is an open or closed dual circuit with a positive probability in any annulus
$A(n, 2n)$. We still call this argument the RSW lemma.

 By using (43)--(53) in Kesten (1986), we know that for any $l\geq 1$ there exists $C=C(l)>0$ such that
$$E_\nu T_n^l \leq C (E_\nu T_n)^l= n^{l(2-5/48+o(1))}.\eqno{(2.1)}$$
Similarly, by the exact proof in Kesten's  above estimate, we still have
$$E_{p_c} S_n^l \leq C (E_{p_c} S_n)^l= n^{l(2-5/48+o(1))}.\eqno{(2.2)}$$

By (2.2) and Cauchy-Schwarz's inequality, we have the followng lemma.\\

{\bf Lemma 2.2.} {\em For any $\epsilon >0$, there exists a positive constant $C_1$ independent of  $n$ such that
$$ P_{p_c} ( S_n \geq n^{2-5/48-\epsilon}) \geq C_1.\eqno{}$$}

Given  four-arm paths in a square $B(n)$ and four-arm paths in an annulus  $A(n, m)$ for $n< m$, one of the major estimates (see Lemma 4 and Lemma 6 in Kesten (1987)) is to reconnect them  by costing a constant probability independent of $n$ and $m$.  It is a more general argument than the RSW lemma since the RSW lemma only reconnects the open or closed  paths. This connection  is called the {\em reconnection lemma}. \\

{\bf Reconnection lemma.} (Kesten (1987)). {\em For $n \leq m$,
there exists $C>0$ (independent of $n,m$) such that}
$$ P_{p_c}( {\cal Q}_4(n) P_p({\cal Q}_4(n,m)))\leq C P_{p_c}({\cal Q}_4(m)).$$
\\

In fact, in Kesten's reconnection lemma, if there are four-arm paths, then  we can restrict the four-arm paths  in  special locations. More precisely,  we denote by ${\cal E}_4(k)$ the event that there are four-arm paths from some edge $e$ and $e^*$  in $B(2^{k-1})$ to  $\partial B(2^k)$.
We define the sub-event 
$\bar{\cal E}_4(2^{k})$ of ${\cal E}_4( 2^k)$  such that there are four-arm paths $r_1, r_2^*,
r_3,$ and $ r_4^*$ from $e$ and $e^*$  in ${\cal E}_4(2^k)$ with (see Fig. 1)
\begin{eqnarray*}
&& r_1\cap A(2^{k-1}, 2^k) \subset {\bf A}_1:=[-2^k, -2^{k-1}]\times [-2^{k-1}, 2^{k-1}],\\
&& r_3 \cap A(2^{k-1}, 2^k) \subset {\bf A}_3:=[2^{k-1}, 2^{k}]\times [-2^{k-1}, 2^{k-1}],\\
&& r_2^* \cap A(2^{k-1}, 2^k) \subset {\bf A}_2:=[-2^{k-1}, 2^{k-1}]\times [2^{k-1}, 2^{k}],\\
&& r_4^* \cap A(2^{k-1}, 2^k) \subset {\bf A}_4:=[-2^{k-1}, 2^{k-1}]\times [-2^{k}, -2^{k-1}].
\end{eqnarray*}
In addition, we require that there exist an open vertical crossing on ${\bf A}_i$ for $i=1,3$ and 
a horizontal vacant dual crossing on ${\bf A}^*_i$ for $i=2,4$.
These extra vertical and horizontal crossings are called {\em wings}.
Kesten's Lemma 4 (see Kesten (1987)) showed that  there exists $C>0$ such that for all $k$,
$$ P_{p_c}( {\cal E}_4( 2^k)) \leq  C P_{p_c}(\bar{\cal E}_4(2^k)).\eqno{(2.3)}$$
Here we want to remark that Kesten's proof in (2.3) is  for four-arm paths from a fixed edge.
However, the exact proof can show our situation.

It is also easier to handle a one arm path in this way.  We may replace ${\cal E}_4(k)$
and $\bar{\cal E}_4(k)$ with  ${\cal E}_1(k)$
and $\bar{\cal E}_1(k)$ by ignoring the paths $r_2^*$, $r_3$, and $r_4^*$. We can use the same argument of  (2.3) to show that
$$ P_{p_c}( {\cal E}_1( 2^k)) \leq  C P_{p_c}(\bar{\cal E}_1(2^k)).\eqno{(2.4)}$$
For $n< m$, we also define  ${\cal G}_1(k)$ to be the event that there is an open path $r_1$
from the origin to $\partial B(2^k)$. In addition,  let $\bar{\cal G}_1(k)$ be the sub-event of
${\cal G}_1(k)$ such that $r_1$ also crosses ${\bf A}_1$ vertically. It  also follows from the same argument
of (2.4) that
$$ P_{p_c}( {\cal G}_1( 2^k)) \leq  C P_{p_c}(\bar{\cal G}_1(2^k)).\eqno{(2.5)}$$
We want to point out that the proof of (2.5) can be shown directly by using  the FKG inequality and the RSW lemma since we only need to deal with open paths, positive events. The proofs of 
Kesten's Lemma 4  is much more  complicated since they need to deal with
both open and closed paths, positive and negative events.

 For each  $e \in B(n)$ with two vertices $v'=v'(e)$ and $v''=v''(e)$,  $e$ is called a {\em pivotal} for the connection of 
 $v'(e)$ and the origin if there is or  isn't an open path inside
 B(n) from $v'(e)$ to the origin  if $e$ is open or closed.  By Lemma 2.1, it is easy to 
 verify that  if $e$ is pivotal, then there is  a  dual path $D^*(v')$ from one vertex of $e^*$ to the other
 inside $B(n)$  such that all of its edges, except $e^*$ and its edges on $\partial B(n)$,  are closed, and there is an open path from
 $v''(e)$ to the origin (see Fig. 2). 
 Together with $e^*$, $D^*(v')\cup e^*$ is a circuit, and the circuit separates $v'(e)$ from $v''(e)$.
Note that there might be many such circuits. We always select one containing the smallest
number of edges. It follows from Lemma 2.1 (see Fig. 2) that  $D(v')$ is the exterior boundary of 
the open cluster of $v'(e)$. 
We denote by  $\bar D(v')$ the open cluster of $v'(e)$.  Let ${\cal H}_\epsilon(e)$ be the event that   
$e$ is pivotal for $v'(e)$,  and $D^*(v')$ does not contain the edges in $\partial B^*(n)$,
 and $| B(n) \cap  \bar{D}(v')|\geq n^{2-5/48-\epsilon}$.   Note that $e$ can be either open or closed independent of ${\cal H}_\epsilon(e)$.  
 On ${\cal H}_\epsilon(e)$, if $e$ is changed from open to closed, then ${\cal C}$ will lose at least $n^{2-5/48-\epsilon}$ many vertices. \\
 
 {\bf Lemma 2.3.}
{\em If ${\cal H}_\epsilon =\cup_{e\in B(n)} {\cal H}(e)$, then  for each $\epsilon >0$
there is a constant $C>0$ (independent  of $\epsilon$ and $n$)  such that
$$ P_\nu ( {\cal H}_\epsilon) \geq C .$$
In addition, there exists a constant $C_i$ for $i=1,2$ such that
$$ 0.5 P_\nu ( {\cal H}_\epsilon(e))= P_\nu ( {\cal H}_\epsilon(e), \mbox{ $e$ is open} ) =P_\nu ( {\cal H}_\epsilon(e), \mbox{ $e$ is closed} ).$$}

{\bf Proof.}  We focus on $S=[n/16, 3n/16]^2$.  Let $v$ be the center of $S$. By using the RSW lemma, there exists an edge $e$ in $v+ [n/8, n/8]^2$ such that 
${\cal E}_4(v, n/8)$ occurs with a positive probability, where ${\cal E}_4(v, n/8)$
is the same event  ${\cal E}_4(n/8)$ by shifting each configuration from the origin to $v$.
That is,
$$ P_{p_c} ( {\cal E}_4(v, n/8)) \geq C_1.\eqno{(2.6)}$$
In fact,  Kesten,  Sidoravicius,  and Zhang (1998) can show that there exist five-arm paths with a positive probability.
By (2.4),  we have
$$ P_{p_c} ( \bar{\cal E}_4(v, n/8)) \geq C_1.\eqno{(2.7)}$$
Here we only need the closed wings in ${\bf A}_2$ and ${\bf A}_4$.
On $ \bar{\cal E}_4(v, n/8)$, we define the following events (see Fig. 1).
Let ${\cal D}_1(n)$ be the event that there is a closed dual path
from $[n/16, 3n/16]\times \{3n/16\}$ to $[n/16, 3n/16]\times \{n/16\}$ inside  annulus
$${\bf B}_1=[n/16, 11n/16]\times [-3n/16, 7n/16]\setminus[3n/16, 9n/16]\times [-n/16, 5n/16].$$
Let  ${\cal D}_2(n)$ be the event that there is an open circuit in annulus 
$${\bf B}_2=[3n/16, 9n/16] \times
[-n/16, 5n/16]\setminus [5n/16, 7n/16]\times [n/16, 3n/16].$$
On ${\cal D}_2(n)$, there might be many such circuits.
We assume that $D_2(n)$ is the largest circuit among all the circuits.
By Proposition 2.3  in Kesten (1982), if $D_2(n)=\Gamma$ for a fixed  circuit $\Gamma$,
then  
$$\mbox{ configurations of edges inside $\dot{\Gamma}$ and outside $\dot{\Gamma}$ are independent,}\eqno{(2.8)}$$
where $\dot{\Gamma}$ is the edges enclosed by $\Gamma$, but not in $\Gamma$.
On $D_2(n)=\Gamma$ for a fixed $\Gamma$, let  ${\cal N}(\Gamma)$ be the event that there are more than 
$n^{2-5/48-\epsilon}$ many vertices inside $[5n/16, 7n/16]\times [n/16, 3n/16]$ such that each is connected to $\Gamma$  by open paths.  Let ${\cal N}(D_2(n))= \cup_{\Gamma} {\cal N}(\Gamma)$. By Lemma 2.2, there exists $C_2>0$ such that
$$P_{p_c}( {\cal N}(\Gamma) \,\, |\,\, D_2(n)=\Gamma) \geq P_{p_c} (S_{n/8} \geq n^{5/48-\epsilon} ) \geq C_2.\eqno{( 2.9)}$$
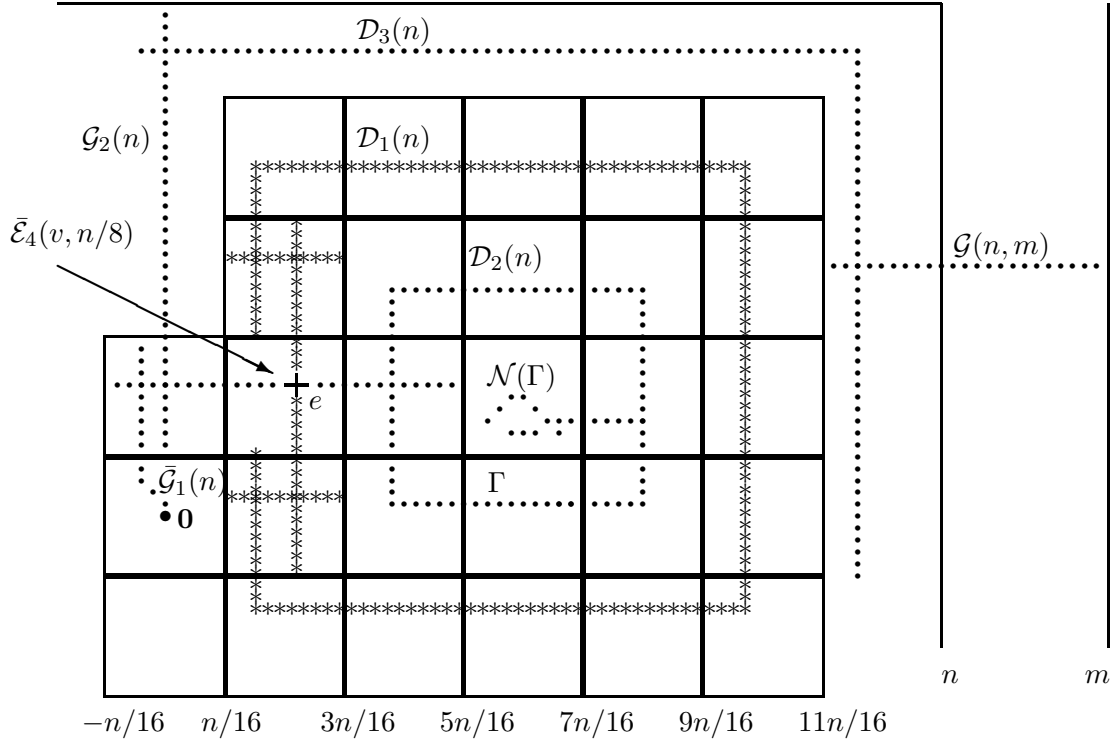
\begin{figure}
\begin{center}
\setlength{\unitlength}{0.0125in}%
\begin{picture}(250,220)(67,770)
\thicklines
\put(0,770){\framebox(50,50)[br]{$$}}
\put(50,770){\framebox(50,50)[br]{$$}}
\put(100,770){\framebox(50,50)[br]{$$}}
\put(150,770){\framebox(50,50)[br]{$$}}
\put(200,770){\framebox(50,50)[br]{$$}}
\put(250,770){\framebox(50,50)[br]{$$}}
\put(0,870){\framebox(50,50)[br]{$$}}
\put(0,820){\framebox(50,50)[br]{$$}}
\put(50,820){\framebox(50,50)[br]{$$}}
\put(50,870){\framebox(50,50)[br]{$$}}
\put(50,920){\framebox(50,50)[br]{$$}}
\put(150,820){\framebox(50,50)[br]{$$}}
\put(150,870){\framebox(50,50)[br]{$$}}
\put(150,920){\framebox(50,50)[br]{$$}}
\put(100,820){\framebox(50,50)[br]{$$}}
\put(100,870){\framebox(50,50)[br]{$$}}
\put(100,920){\framebox(50,50)[br]{$$}}
\put(250,820){\framebox(50,50)[br]{$$}}
\put(250,870){\framebox(50,50)[br]{$$}}
\put(250,920){\framebox(50,50)[br]{$$}}
\put(200,820){\framebox(50,50)[br]{$$}}
\put(200,870){\framebox(50,50)[br]{$$}}
\put(200,920){\framebox(50,50)[br]{$$}}
\put(250,820){\framebox(50,50)[br]{$$}}
\put(50,970){\framebox(50,50)[br]{$$}}
\put(100,970){\framebox(50,50)[br]{$$}}
\put(200,970){\framebox(50,50)[br]{$$}}
\put(150,970){\framebox(50,50)[br]{$$}}
\put(250,970){\framebox(50,50)[br]{$$}}
\put(-20,1060){\line(1,0){370}}
\put(350,1060){\line(0, -1){270}}
\put(420,1060){\line(0, -1){270}}

\put(75,900){\line(1,0){10}}
\put(80,895){\line(0,1){10}}
\put(85,890){\mbox{$e$}{}}
\put(-40,960){\mbox{$\bar{\cal E}_4(v,n/8)$}{}}
\put(-20,950){\vector(2,-1){90}}
\put(105,1000){\mbox{${\cal D}_1(n)$}{}}
\put(105,1045){\mbox{${\cal D}_3(n)$}{}}
\put(152,950){\mbox{${\cal D}_2(n)$}{}}
\put(22,856){\mbox{$\bar{\cal G}_1(n)$}{}}
\put(-10,1000){\mbox{${\cal G}_2(n)$}{}}
\put(-10,755){\mbox{$-n/16$}{}}
\put(40,755){\mbox{$n/16$}{}}
\put(90,755){\mbox{$3n/16$}{}}
\put(140,755){\mbox{$5n/16$}{}}
\put(190,755){\mbox{$7n/16$}{}}
\put(240,755){\mbox{$9n/16$}{}}
\put(290,755){\mbox{$11n/16$}{}}
\put(350,775){\mbox{$n$}{}}
\put(410,775){\mbox{$m$}{}}
\put(355,955){\mbox{${\cal G}(n,m)$}{}}
\put(160,855){\mbox{$\Gamma$}{}}

\put(90,900){\circle*{2}}
\put(95,900){\circle*{2}}
\put(100,900){\circle*{2}}
\put(105,900){\circle*{2}}
\put(110,900){\circle*{2}}
\put(115,900){\circle*{2}}
\put(120,900){\circle*{2}}
\put(125,900){\circle*{2}}
\put(130,900){\circle*{2}}
\put(135,900){\circle*{2}}
\put(140,900){\circle*{2}}
\put(145,900){\circle*{2}}

\put(70,900){\circle*{2}}
\put(65,900){\circle*{2}}
\put(60,900){\circle*{2}}
\put(55,900){\circle*{2}}
\put(50,900){\circle*{2}}
\put(45,900){\circle*{2}}
\put(40,900){\circle*{2}}
\put(35,900){\circle*{2}}
\put(30,900){\circle*{2}}
\put(25,900){\circle*{2}}
\put(20,900){\circle*{2}}
\put(15,900){\circle*{2}}
\put(10,900){\circle*{2}}
\put(5,900){\circle*{2}}

\put(120,905){\circle*{2}}
\put(120,910){\circle*{2}}
\put(120,915){\circle*{2}}
\put(120,920){\circle*{2}}
\put(120,925){\circle*{2}}
\put(120,930){\circle*{2}}
\put(120,935){\circle*{2}}
\put(120,940){\circle*{2}}
\put(125,940){\circle*{2}}
\put(130,940){\circle*{2}}
\put(135,940){\circle*{2}}
\put(140,940){\circle*{2}}
\put(145,940){\circle*{2}}
\put(150,940){\circle*{2}}
\put(155,940){\circle*{2}}
\put(160,940){\circle*{2}}
\put(165,940){\circle*{2}}
\put(170,940){\circle*{2}}
\put(175,940){\circle*{2}}
\put(180,940){\circle*{2}}
\put(185,940){\circle*{2}}
\put(190,940){\circle*{2}}
\put(195,940){\circle*{2}}
\put(200,940){\circle*{2}}
\put(205,940){\circle*{2}}
\put(210,940){\circle*{2}}
\put(215,940){\circle*{2}}
\put(220,940){\circle*{2}}
\put(225,940){\circle*{2}}
\put(225,935){\circle*{2}}
\put(225,930){\circle*{2}}
\put(225,925){\circle*{2}}
\put(225,920){\circle*{2}}
\put(225,915){\circle*{2}}
\put(225,910){\circle*{2}}
\put(225,905){\circle*{2}}
\put(225,900){\circle*{2}}
\put(225,895){\circle*{2}}
\put(225,890){\circle*{2}}
\put(225,885){\circle*{2}}
\put(225,880){\circle*{2}}
\put(225,875){\circle*{2}}
\put(225,870){\circle*{2}}
\put(225,865){\circle*{2}}
\put(225,860){\circle*{2}}
\put(225,855){\circle*{2}}
\put(225,850){\circle*{2}}
\put(220,850){\circle*{2}}
\put(215,850){\circle*{2}}
\put(210,850){\circle*{2}}
\put(205,850){\circle*{2}}
\put(200,850){\circle*{2}}
\put(195,850){\circle*{2}}
\put(190,850){\circle*{2}}
\put(195,850){\circle*{2}}
\put(190,850){\circle*{2}}
\put(185,850){\circle*{2}}
\put(180,850){\circle*{2}}
\put(175,850){\circle*{2}}
\put(170,850){\circle*{2}}
\put(165,850){\circle*{2}}
\put(160,850){\circle*{2}}
\put(155,850){\circle*{2}}
\put(150,850){\circle*{2}}
\put(145,850){\circle*{2}}
\put(140,850){\circle*{2}}
\put(135,850){\circle*{2}}
\put(130,850){\circle*{2}}
\put(125,850){\circle*{2}}
\put(120,850){\circle*{2}}
\put(120,855){\circle*{2}}
\put(120,860){\circle*{2}}
\put(120,865){\circle*{2}}
\put(120,870){\circle*{2}}
\put(120,875){\circle*{2}}
\put(120,880){\circle*{2}}
\put(120,885){\circle*{2}}
\put(120,890){\circle*{2}}
\put(120,895){\circle*{2}}

\put(170,895){\circle*{2}}
\put(165,890){\circle*{2}}
\put(160,885){\circle*{2}}
\put(170,880){\circle*{2}}
\put(175,880){\circle*{2}}
\put(175,895){\circle*{2}}
\put(180,890){\circle*{2}}
\put(185,885){\circle*{2}}
\put(180,880){\circle*{2}}
\put(190,880){\circle*{2}}
\put(190,885){\circle*{2}}
\put(195,885){\circle*{2}}
\put(200,885){\circle*{2}}
\put(205,885){\circle*{2}}
\put(210,885){\circle*{2}}
\put(215,885){\circle*{2}}
\put(220,885){\circle*{2}}
\put(225,885){\circle*{2}}
\put(160,900){\mbox{${\cal N}(\Gamma)$}}
\put(25,850){\circle*{2}}
\put(20,855){\circle*{2}}
\put(15,860){\circle*{2}}
\put(15,865){\circle*{2}}
\put(15,870){\circle*{2}}
\put(15,875){\circle*{2}}
\put(15,880){\circle*{2}}
\put(15,885){\circle*{2}}
\put(15,890){\circle*{2}}
\put(15,895){\circle*{2}}
\put(15,900){\circle*{2}}
\put(15,905){\circle*{2}}
\put(15,910){\circle*{2}}
\put(15,915){\circle*{2}}

\put(25,870){\circle*{2}}
\put(25,875){\circle*{2}}
\put(25,880){\circle*{2}}
\put(25,885){\circle*{2}}
\put(25,890){\circle*{2}}
\put(25,895){\circle*{2}}
\put(25,900){\circle*{2}}
\put(25,905){\circle*{2}}
\put(25,910){\circle*{2}}
\put(25,915){\circle*{2}}
\put(25,920){\circle*{2}}
\put(25,925){\circle*{2}}
\put(25,930){\circle*{2}}
\put(25,935){\circle*{2}}
\put(25,940){\circle*{2}}
\put(25,945){\circle*{2}}
\put(25,950){\circle*{2}}
\put(25,955){\circle*{2}}
\put(25,960){\circle*{2}}
\put(25,965){\circle*{2}}
\put(25,970){\circle*{2}}
\put(25,975){\circle*{2}}
\put(25,980){\circle*{2}}
\put(25,985){\circle*{2}}
\put(25,990){\circle*{2}}
\put(25,995){\circle*{2}}
\put(25,1000){\circle*{2}}
\put(25,1005){\circle*{2}}
\put(25,1010){\circle*{2}}
\put(25,1015){\circle*{2}}
\put(25,1020){\circle*{2}}
\put(25,1025){\circle*{2}}
\put(25,1030){\circle*{2}}
\put(25,1035){\circle*{2}}
\put(25,1040){\circle*{2}}
\put(25,1045){\circle*{2}}
\put(25,1050){\circle*{2}}
\put(25,1055){\circle*{2}}
\put(30,1040){\circle*{2}}
\put(35,1040){\circle*{2}}
\put(20,1040){\circle*{2}}
\put(15,1040){\circle*{2}}
\put(40,1040){\circle*{2}}
\put(45,1040){\circle*{2}}
\put(50,1040){\circle*{2}}
\put(55,1040){\circle*{2}}
\put(60,1040){\circle*{2}}
\put(65,1040){\circle*{2}}
\put(70,1040){\circle*{2}}
\put(75,1040){\circle*{2}}
\put(80,1040){\circle*{2}}
\put(85,1040){\circle*{2}}
\put(90,1040){\circle*{2}}
\put(95,1040){\circle*{2}}
\put(100,1040){\circle*{2}}
\put(105,1040){\circle*{2}}
\put(110,1040){\circle*{2}}
\put(115,1040){\circle*{2}}
\put(120,1040){\circle*{2}}
\put(125,1040){\circle*{2}}
\put(130,1040){\circle*{2}}
\put(135,1040){\circle*{2}}
\put(140,1040){\circle*{2}}
\put(145,1040){\circle*{2}}
\put(150,1040){\circle*{2}}
\put(155,1040){\circle*{2}}
\put(160,1040){\circle*{2}}
\put(165,1040){\circle*{2}}
\put(170,1040){\circle*{2}}
\put(175,1040){\circle*{2}}
\put(180,1040){\circle*{2}}
\put(185,1040){\circle*{2}}
\put(190,1040){\circle*{2}}
\put(195,1040){\circle*{2}}
\put(200,1040){\circle*{2}}
\put(205,1040){\circle*{2}}
\put(210,1040){\circle*{2}}
\put(215,1040){\circle*{2}}
\put(220,1040){\circle*{2}}
\put(225,1040){\circle*{2}}
\put(230,1040){\circle*{2}}
\put(235,1040){\circle*{2}}
\put(240,1040){\circle*{2}}
\put(245,1040){\circle*{2}}
\put(250,1040){\circle*{2}}
\put(255,1040){\circle*{2}}
\put(260,1040){\circle*{2}}
\put(265,1040){\circle*{2}}
\put(270,1040){\circle*{2}}
\put(275,1040){\circle*{2}}
\put(280,1040){\circle*{2}}
\put(285,1040){\circle*{2}}
\put(290,1040){\circle*{2}}
\put(295,1040){\circle*{2}}
\put(300,1040){\circle*{2}}
\put(305,1040){\circle*{2}}
\put(310,1040){\circle*{2}}
\put(315,1040){\circle*{2}}
\put(315,1035){\circle*{2}}
\put(315,1030){\circle*{2}}
\put(315,1025){\circle*{2}}
\put(315,1020){\circle*{2}}
\put(315,1015){\circle*{2}}
\put(315,1010){\circle*{2}}
\put(315,1005){\circle*{2}}
\put(315,1000){\circle*{2}}
\put(315,995){\circle*{2}}
\put(315,990){\circle*{2}}
\put(315,985){\circle*{2}}
\put(315,980){\circle*{2}}
\put(315,975){\circle*{2}}
\put(315,970){\circle*{2}}
\put(315,965){\circle*{2}}
\put(315,960){\circle*{2}}
\put(315,955){\circle*{2}}
\put(315,950){\circle*{2}}
\put(315,945){\circle*{2}}
\put(315,940){\circle*{2}}
\put(315,935){\circle*{2}}
\put(315,930){\circle*{2}}
\put(315,925){\circle*{2}}
\put(315,920){\circle*{2}}
\put(315,915){\circle*{2}}
\put(315,910){\circle*{2}}
\put(315,905){\circle*{2}}
\put(315,900){\circle*{2}}
\put(315,895){\circle*{2}}
\put(315,890){\circle*{2}}
\put(315,885){\circle*{2}}
\put(315,880){\circle*{2}}
\put(315,875){\circle*{2}}
\put(315,870){\circle*{2}}
\put(315,865){\circle*{2}}
\put(315,860){\circle*{2}}
\put(315,855){\circle*{2}}
\put(315,850){\circle*{2}}
\put(315,845){\circle*{2}}
\put(315,840){\circle*{2}}
\put(315,835){\circle*{2}}
\put(315,830){\circle*{2}}
\put(315,825){\circle*{2}}
\put(315,820){\circle*{2}}

\put(305,950){\circle*{2}}
\put(310,950){\circle*{2}}
\put(315,950){\circle*{2}}
\put(320,950){\circle*{2}}
\put(325,950){\circle*{2}}
\put(330,950){\circle*{2}}
\put(335,950){\circle*{2}}
\put(340,950){\circle*{2}}
\put(345,950){\circle*{2}}
\put(350,950){\circle*{2}}
\put(355,950){\circle*{2}}
\put(360,950){\circle*{2}}
\put(365,950){\circle*{2}}
\put(370,950){\circle*{2}}
\put(375,950){\circle*{2}}
\put(380,950){\circle*{2}}
\put(385,950){\circle*{2}}
\put(390,950){\circle*{2}}
\put(395,950){\circle*{2}}
\put(400,950){\circle*{2}}
\put(405,950){\circle*{2}}
\put(410,950){\circle*{2}}
\put(415,950){\circle*{2}}

\put(77,902){*}
\put(77,907){*}
\put(77,912){*}
\put(77,917){*}
\put(77,922){*}
\put(77,927){*}
\put(77,932){*}
\put(77,937){*}
\put(77,942){*}
\put(77,947){*}
\put(77,952){*}
\put(77,957){*}
\put(77,962){*}

\put(77,887){*}
\put(77,882){*}
\put(77,877){*}
\put(77,872){*}
\put(77,867){*}
\put(77,862){*}
\put(77,857){*}
\put(77,852){*}
\put(77,847){*}
\put(77,842){*}
\put(77,837){*}
\put(77,832){*}
\put(77,827){*}
\put(77,822){*}
\put(77,817){*}

\put(50,847){*}
\put(55,847){*}
\put(60,847){*}
\put(65,847){*}
\put(70,847){*}
\put(75,847){*}
\put(80,847){*}
\put(85,847){*}
\put(90,847){*}
\put(95,847){*}

\put(50,947){*}
\put(55,947){*}
\put(60,947){*}
\put(65,947){*}
\put(70,947){*}
\put(75,947){*}
\put(80,947){*}
\put(85,947){*}
\put(90,947){*}
\put(95,947){*}

\put(60,915){*}
\put(60,920){*}
\put(60,925){*}
\put(60,930){*}
\put(60,935){*}
\put(60,940){*}
\put(60,945){*}
\put(60,950){*}
\put(60,955){*}
\put(60,960){*}
\put(60,965){*}
\put(60,970){*}
\put(60,975){*}
\put(60,980){*}
\put(60,985){*}
\put(65,985){*}
\put(70,985){*}
\put(75,985){*}
\put(80,985){*}
\put(85,985){*}
\put(90,985){*}
\put(95,985){*}
\put(100,985){*}
\put(105,985){*}
\put(110,985){*}
\put(115,985){*}
\put(120,985){*}
\put(125,985){*}
\put(130,985){*}
\put(135,985){*}
\put(140,985){*}
\put(145,985){*}
\put(150,985){*}
\put(155,985){*}
\put(160,985){*}
\put(165,985){*}
\put(170,985){*}
\put(175,985){*}
\put(180,985){*}
\put(185,985){*}
\put(190,985){*}
\put(195,985){*}
\put(200,985){*}
\put(205,985){*}
\put(210,985){*}
\put(215,985){*}
\put(220,985){*}
\put(225,985){*}
\put(230,985){*}
\put(235,985){*}
\put(240,985){*}
\put(245,985){*}
\put(250,985){*}
\put(255,985){*}
\put(260,985){*}
\put(265,985){*}
\put(265,980){*}
\put(265,975){*}
\put(265,970){*}
\put(265,965){*}
\put(265,960){*}
\put(265,955){*}
\put(265,950){*}
\put(265,945){*}
\put(265,940){*}
\put(265,935){*}
\put(265,930){*}
\put(265,925){*}
\put(265,920){*}
\put(265,915){*}
\put(265,910){*}
\put(265,905){*}
\put(265,900){*}
\put(265,895){*}
\put(265,890){*}
\put(265,885){*}
\put(265,880){*}
\put(265,875){*}
\put(265,870){*}
\put(265,865){*}
\put(265,860){*}
\put(265,855){*}
\put(265,850){*}
\put(265,845){*}
\put(265,840){*}
\put(265,835){*}
\put(265,830){*}
\put(265,825){*}
\put(265,820){*}
\put(265,815){*}
\put(265,810){*}
\put(265,805){*}
\put(265,800){*}
\put(260,800){*}
\put(255,800){*}
\put(250,800){*}
\put(245,800){*}
\put(240,800){*}
\put(235,800){*}
\put(230,800){*}
\put(225,800){*}
\put(220,800){*}
\put(215,800){*}
\put(210,800){*}
\put(205,800){*}
\put(200,800){*}
\put(195,800){*}
\put(190,800){*}
\put(185,800){*}
\put(180,800){*}
\put(175,800){*}
\put(170,800){*}
\put(165,800){*}
\put(160,800){*}
\put(155,800){*}
\put(150,800){*}
\put(145,800){*}
\put(140,800){*}
\put(135,800){*}
\put(130,800){*}
\put(125,800){*}
\put(120,800){*}
\put(115,800){*}
\put(110,800){*}
\put(105,800){*}
\put(100,800){*}
\put(95,800){*}
\put(90,800){*}
\put(85,800){*}
\put(80,800){*}
\put(75,800){*}
\put(70,800){*}
\put(65,800){*}
\put(60,800){*}
\put(60,805){*}
\put(60,810){*}
\put(60,815){*}
\put(60,820){*}
\put(60,825){*}
\put(60,830){*}
\put(60,835){*}
\put(60,840){*}
\put(60,845){*}
\put(60,850){*}
\put(60,855){*}
\put(60,860){*}
\put(60,865){*}

\put(25,845){\circle*{4}}
\put(30,840){\mbox{${\bf 0}$}{}}
\thicklines
\end{picture}
\end{center}
\caption{\em  The figure shows how to construct a few events together such that
${\cal H}_\epsilon$ occurs. ${\cal N}(\Gamma)$ is all vertices in $[5n/16, 7n/16]\times [n/16, 3n/16]$ such that each of them is connected to $\Gamma$ by an open path. ${\cal E}_4(v, n/8)$ is the event  that there are four-arm paths from $e$ and these closed paths have two wings.}
\end{figure}

Let ${\cal G}_2(n)$ be the event that there is an  open vertical crossing   in $[-n/16, n/16]\times 
[n/16, n]$.
Let ${\cal D}_3(n)$ be the event that there is an open circuit in the annulus $A(11n/16, n)$, and let
${\cal G}(n, m) $  be the event that  there is an open path  from $ B(11n/16)$ to $\partial B(m)$.
Now we put these events together.  For any $\delta >0$, we take $m$ large such that 
\begin{eqnarray*}
&&P_\nu( {\cal H}_\epsilon) = \lim_{m\rightarrow \infty} P_{p_c}^{-1} ({\bf 0} \rightarrow \partial B(m))P_{p_c} ({\cal H}_\epsilon \cap {\bf 0} \rightarrow \partial B(m))\\
&\geq & P_{p_c}^{-1} ({\bf 0} \rightarrow \partial B(m))P_{p_c} ({\cal H}_\epsilon \cap {\bf 0} \rightarrow \partial B(m))-\delta.  \hskip 8cm  (2.10)
\end{eqnarray*}
Let us focus on $P_{p_c} ({\cal H}_\epsilon \cap {\bf 0} \rightarrow \partial B(m))$. 
As we defined these events (see Fig. 1),
\begin{eqnarray*}
&&P_{p_c} ({\cal H}_\epsilon \cap {\bf 0} \rightarrow \partial B(m))\\
&\geq  &P_{p_c} (\bar{\cal E}_4(n/8) \cap {\cal D}_1(n) \cap {\cal D}_2(n) \cap {\cal D}_3(n) \cap {\cal N}(D_2(n))\cap \bar{\cal G}_1(n)\cap {\cal G}_2(n)\cap{\cal G}(n, m) ).\hskip 2.5cm (2.11)
\end{eqnarray*}
Note that  open paths and closed paths are positive and negative events, so by Lemma 3 in Kesten (1987),
\begin{eqnarray*}
&&  P_{p_c} (\bar{\cal E}_4(n/8) \cap {\cal D}_1(n) \cap {\cal D}_2(n) \cap {\cal D}_3(n) \cap {\cal N}(D_2(n))\cap \bar{\cal G}_1(n)\cap {\cal G}_2(n)\cap{\cal G}(n, m) )\\
&\geq & P_{p_c}( \bar{\cal E}_4(n/8) )P_{p_c}(  {\cal D}_1(n) \cap {\cal D}_2(n)  \cap{\cal D}_3(n) \cap {\cal N}(D_n(n))\cap \bar{\cal G}_1(n)\cap {\cal G}_2(n)\cap{\cal G}(n, m) ).
\hskip 2cm (2.12)
\end{eqnarray*}
Note that ${\cal D}_1$ and ${\cal D}_2\cap {\cal D}_3\cap {\cal N}(D_2)\cap \bar{\cal G}_1(n)$
occur in different areas (see Fig. 1), so by (2.7) together with   the RSW lemma,
\begin{eqnarray*}
&&P_{p_c}(  {\cal D}_1(n) \cap {\cal D}_2(n)  \cap {\cal D}_3(n) \cap {\cal N}(D_n(n))\cap \bar{\cal G}_1(n)\cap {\cal G}_2(n)\cap{\cal G}(n, m) )\\
&\geq & C P_{p_c} ({\cal D}_2(n)  \cap {\cal D}_3(n) \cap {\cal N}(D_n(n))\cap \bar{\cal G}_1(n)\cap {\cal G}_2(n)\cap{\cal G}(n, m) ).\hskip 5cm {(2.13)}
\end{eqnarray*}
By (2.9),
\begin{eqnarray*}
&&P_{p_c} ({\cal D}_2(n) \cap  {\cal D}_3(n) \cap {\cal N}(D_n(n))\cap \bar{\cal G}_1(n)\cap {\cal G}_2(n)\cap{\cal G}(n, m) )\\
&=& \sum_{\Gamma} P_{p_c} ( {\cal N}(\Gamma) \,\, |\,\, D_2(n)=\Gamma)  P_{p_c} ( {\cal D}_3(n)\cap \bar{\cal G}_1(n)\cap {\cal G}_2(n)\cap{\cal G}(n, m) \cap \{D_2(n)=\Gamma\})\\
&\geq & C_2 \sum_{\Gamma}  P_{p_c} ( {\cal D}_3(n)\cap \bar{\cal G}_1(n)\cap{\cal G}_2\cap{\cal G}(n, m) \cap \{D_2(n)=\Gamma\})\\
&\geq & C_2  P_{p_c} ( {\cal D}_3(n)\cap \bar{\cal G}_1(n)\cap {\cal G}_2(n)\cap{\cal G}(n, m) \cap{\cal D}_2(n)),\hskip 7cm (2.14)
\end{eqnarray*}
where the sum takes over all possible edge sets $\Gamma$.
By the FKG inequality,  the RSW lemma, and  (2.5),
$$P_{p_c} ( {\cal D}_3(n)\cap\bar{\cal G}_1(n)\cap{\cal G}_2(n)\cap{\cal G}(n, m) \cap{\cal D}_2(n))\geq C_3 P_{p_c} ( {\cal G}_1(n)\cap{\cal G}(n, m))\geq 
C_3 P_{p_c}({\bf 0} \rightarrow \partial B(m)).\eqno{(2.15)}$$
If we put (2.11)--(2.15) together, then 
$$P_\nu({\cal H}_\epsilon) \geq C_3-\delta.\eqno{(2.16)}$$
Since $\delta $ can be arbitrarily small, the first inequality  in Lemma 2.3 follows. 

Let us show the second inequities in Lemma 2.3.  By using Lemma 2.11, we can show that 
$$P_\nu( {\cal H}_\epsilon(e), e\mbox{ is  open/closed}) = \lim_{m\rightarrow \infty} P_{p_c}^{-1} ({\bf 0} \rightarrow \partial B(m))P_{p_c} ({\cal H}_\epsilon(e), e\mbox{ is  open/closed} , {\bf 0} \rightarrow \partial B(m)).\eqno{(2.17)}$$
Note that the open path from the origin to $\partial B(m)$ will not use $e$, so the configuration of $e$  is open or closed  independent  of $\{{\cal H}_\epsilon(e)\cap {\bf 0} \rightarrow \partial B(m)\}$. Thus,  the second inequality in Lemma 2.3 follows from (2.17).
$\Box$\\

\section{  A martingale construction for $T_n$.}
The martingale construction has been used to estimate  the variance  of an open cluster in  super-critical  percolation (see Zhang 2001). As we know,  super-critical percolation does not have two large open clusters, so it can use the edge-to-edge martingale construction. However, there might be two large open clusters at the criticality. Thus, we need a much more complicated martingale construction. 
We order the edges of $B(n)$  in the following way. First, we order the edges in $ B(1)$ in a certain order. For example, we can do it by  spiral  order.
  Suppose that  we order the edges in $B(i)$; then we can order the edges $B(i+1)\setminus B(i)$ in a certain order.
Thus, we order  all edges in $B(n)$ by $\{e_1, e_2, \cdots, e_i, e_{i+1}, \cdots, e_m\}$.  Our probability space is 
$ \Omega=\{0, 1\}^{m}$, where $m$ is the number of edges in $B(n)$
and,  the generic point of $\Omega$ is $\{\omega_1, \cdots ,\omega_m\}$. In configuration $\omega$, the state of $e_k$ is
$\omega(e_k)=\omega_k$, where  random  variable $\omega_k$ takes value zero or one  if the edge is open or closed.  
Note that $T_n$ is defined on the configuration of open and closed edges in $[-n,n]^2$. If we change the open edge by zero and the closed edge by one, then we can define $T_n(\omega)$ for $\omega\in \Omega$. In particular, if the edges in $B(n)$ are divided into two sub-edge sets $B_1\cup B_2$, then we can write 
$\omega=(\omega(B_1), \omega(B_2))$ and
$T_n(\omega)= T_n(\omega(B_1), \omega(B_2)).$
 
For $e_k$, recall that we considered its two vertices $v_k'$  and $v_k''$ in the last section. 
We also 
let $D^*(v_k')$ be the exterior boundary of the open cluster of $v'$ inside $B(n)$ separating $v_k'$ from $v_k''$ (see Fig. 2).   
  We call $e_k$ the pivotal  in the  connection 
of $v_k'$ and the origin. Now we construct  random sets for each $e_k$. If $e_k$ is a pivotal edge, then 
let $W_k$ be the  edge set enclosed by  $D(v_k')\cup e_k$, including the boundary edges. We also call $e_k$ the pivotal edge of $W_k$. Note that $e_k$ can be open or closed.  If $e_k$ is open or closed, we denote  in particular  by $W_k^+$ or by
$W^-_k$ for $W_k$. If $e_k$ is not a pivotal edge, we simply denote  by $W_k=e_k$.
Note that if both $W_i$ and $W_j$ are non-single edge sets, then the open clusters inside $W_i$ and $W_j$ cannot have a common vertex.  We call this the
{\em disjoint property}.
  We start $e_1$ to construct $W_1$. Suppose that
$\{W_{j_1}, \cdots, W_{j_t}\}$ is defined  by  following the edge
order: $\{e_1, \cdots, e_m\}$. 
Suppose that these edges are 
$$\{e_{i_1},  e_{i_1+1},\cdots, e_{i_1+l_1}, e_{i_2}, e_{i_2+1},  \cdots, e_{i_2+l_2},  \cdots, e_{i_k}, e_{i_k+1}, \cdots, e_{i_k+l_k}\}$$
for $i_j+l_j< i_{j+1}$.
Then we consider $e_{i_1+l_1+1}$ and the random set $W_{i_1+l_1+1}$. We add this random set $W_{i_1+l_1+1}$ to have
$\{W_{j_1}, \cdots, W_{j_t}, W_{j_{t+1}}\}$.
 If we continue this way  by following the edges
$\{e_1, \cdots, e_m\}$, we will construct  all edge sets in $B(n)$.
Without loss of generality, we  still denote them by  $\{W_1, W_2, \cdots, W_m\}$ (see Fig. 2). 
Note that if $W_k$ is not a single edge with the pivotal  edge $e$, then   $D^*(v')$  exists as the exterior boundary 
of the open cluster, denoted by $\bar{D}(v')$ of $v'(e)$. 
We say random set $W_k$ satisfies the {\em bubble condition} if it  indeed consists of the edge set enclosed by $D^*(v_k')\cup e_k^*$ together with $D(v_k')\cup e_k$.
For given sets $\{V_1, \cdots, V_t\}$, 
we may consider the event $\{W_1=V_1, \cdots, W_t=V_t\}$. It follows from our construction that
$$\{W_1=V_1, \cdots, W_t=V_t\}\mbox{ and } \{W_1=V_1', \cdots, W_t=V_t'\} \mbox{ are disjoint if $V_j\neq V_j'$
for some $j\leq t$}. \eqno{(3.1)}$$
In other words,  this construction is unique for each configuration of open and closed edges in $B(n)$. We call this set structure  the {\em spiral  bubble} structure.\\

Note that $W_k$ is a random set, so we use $\omega(W_k)$ to represent the configuration of
$\Omega$  such that  each open or closed edge  in $W_k$ is changed to be zero or one.
 Now we will construct the following martingale sequence
with these sets $\{W_1, W_2, \cdots, W_m \}$. Let filtration 
\begin{figure}
\begin{center}
\setlength{\unitlength}{0.0125in}%
\begin{picture}(250,170)(67,760)
\thicklines

\put(50,770){\framebox(50,50)[br]{$$}}

\put(50,770){\framebox(50,50)[br]{$$}}
\put(100,770){\framebox(50,50)[br]{$$}}
\put(150,770){\framebox(50,50)[br]{$$}}

\put(50,820){\framebox(50,50)[br]{$$}}

\put(150,820){\framebox(50,50)[br]{$$}}
\put(150,870){\framebox(50,50)[br]{$$}}
\put(100,820){\framebox(50,50)[br]{$$}}
\put(100,870){\framebox(50,50)[br]{$$}}
\put(175,795){\framebox(50,50)[br]{$$}}
\put(125,795){\framebox(50,50)[br]{$$}}
\put(225,795){\framebox(50,50)[br]{$$}}
\put(275,795){\framebox(50,50)[br]{$$}}
\put(325,795){\framebox(50,50)[br]{$$}}
\put(325,845){\framebox(50,50)[br]{$$}}
\put(275,845){\framebox(50,50)[br]{$$}}
\put(-20,1060){\line(1,0){395}}
\put(375,1060){\line(0, -1){280}}
\put(139,850){\mbox{$e_i$}{}}
\put(141,828){{\Huge\bf*}}
\put(191,833){{\Huge\bf*}}
\put(241,833){{\Huge\bf*}}
\put(266,858){{\Huge\bf*}}
\put(150,875){\mbox{$v_i''$}{}}
\put(150,826){\mbox{$v_i'$}{}}

\put(116,808){\Huge\bf*}

\put(175,820){\circle*{20}}
\put(225,820){\circle*{20}}
\put(275,820){\circle*{20}}
\put(325,820){\circle*{20}}
\put(175,820){\circle*{20}}
\put(350,845){\circle*{20}}
\put(300,845){\circle*{20}}
\put(325,870){\circle*{20}}
\put(175,820){\line(1,0){200}}
\put(300,820){\framebox(50,50)[br]{$$}}

\put(141,783){{\Huge\bf*}}
\put(191,783){{\Huge\bf*}}
\put(241,783){{\Huge\bf*}}
\put(291,783){{\Huge\bf*}}
\put(341,783){{\Huge\bf*}}
\put(341,883){{\Huge\bf*}}
\put(291,883){{\Huge\bf*}}
\put(90,980){\mbox{$W_j^+$}{}}
\put(50,1000){\mbox{$D^*(v_j')$}{}}
\put(350,1000){\mbox{$B(n)$}{}}

\put(110,915){\mbox{$e_j$}{}}
\put(125,870){\circle*{20}}
\put(100,845){\circle*{20}}
\put(150,895){\circle*{20}}
\put(100,945){\circle*{20}}
\put(125,925){\circle*{20}}
\put(66,908){{\Huge\bf*}}
\put(66,958){{\Huge\bf*}}
\put(91,983){{\Huge\bf*}}
\put(118,958){{\Huge\bf*}}
\put(91,883){{\Huge\bf*}}
\put(355,875){\mbox{$W_i^-$}{}}
\put(310,910){\mbox{$D^*(v_i')$}{}}
\put(75,895){\framebox(50,50)[br]{$$}}
\put(75,945){\framebox(50,50)[br]{$$}}
\put(140,930){\mbox{$v_j''$}{}}
\put(85,925){\mbox{$v_j'$}{}}
\put(100,900){\line(0,1){70}}
\put(90,810){\mbox{${\bf 0}$}{}}
\thicklines
\end{picture}
\end{center}
\caption{\em  The figure shows  how to  construct random sets  $\{W_l\}$ by following the order 
$\{e_0, \cdots, e_m\}$.  $W_1, W_2, \cdots, W_{i-1}$ are single edge sets,  until we reach  the first edge $e_i$ with the vertices $v_i'$ and $v_i''$ such that  $v_i'$ is  enclosed by a closed  dual path $D^*(v_i')\cup e_i^*$, indicted by $*$-edges, and $v_i''$ is connected to the origin by an open path. Note that $D^*(v_i')$ (together with two edges in the boundary) is a  dual circuit separating $v_i'$ from $v_i''$. 
  We add these  edges enclosed by $D^*(v_i')\cup e_i^*$ together with its boundary  edges
to  have $W_i$.   In this case, we denote it by $W^-_i$ since $e_i$ is closed.
We continue to  have  single edge sets until reaching   $e_j$  for $j >i$ such that
 edge $e_j$ also satisfies  the above condition, but with open $e_j$. We then add these edges  enclosed by $D^*(v_j')\cup e_j^*$ together with their boundary  edges to have $W_j$.  
 We denote it by $W_j^+$.
 We continue this  way for all the edges in $B(n)$.
 All the solid-circled edges are open, and all the $*$-dual edges are closed.}
\end{figure}
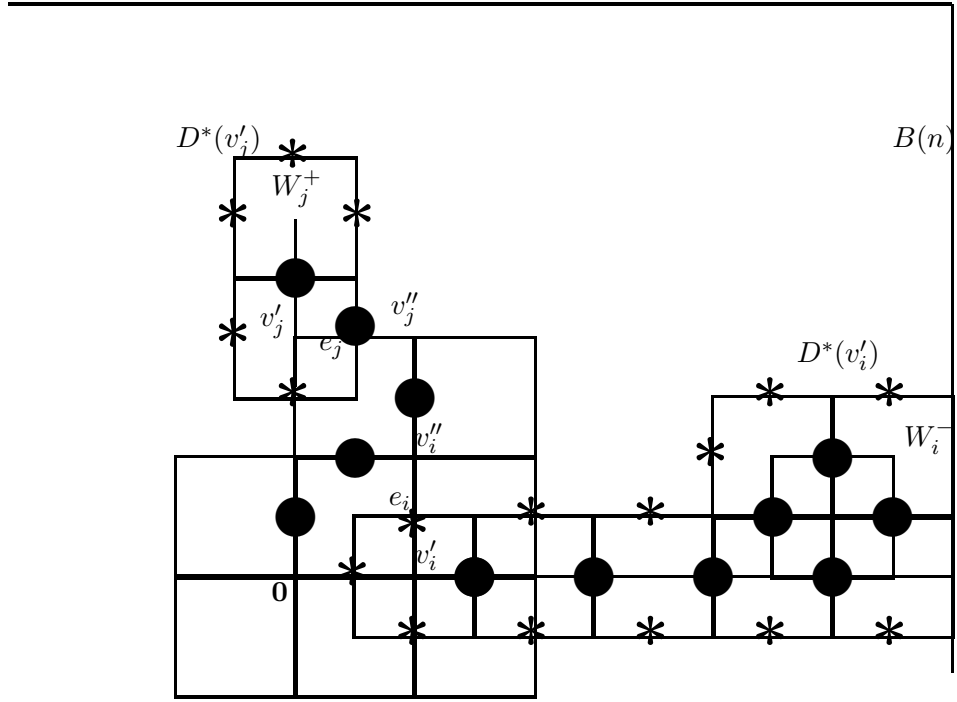

\begin{eqnarray*}
&& {\cal F}_0= \mbox{ the trivial $\sigma$-field }=\{\emptyset, \Omega\}.\\
&& {\cal F}_k= \mbox{  the $\sigma$-field  generated by }\{\omega(e): e\in \bigcup_{i\leq k} W_i\}
\mbox{ for } 0\leq k\leq m.
\end{eqnarray*}
More precisely,  ${\cal F}_k$ is the smallest $\sigma$-field that contains all
the sets with the form 
$$ \{\omega \in \Omega:  \bigcup_{i\leq k} W_i= \{e_{j_1}, e_{j_2}, \cdots, e_{j_t}\}\mbox{ for }  (\omega(e_{j_1}),
\cdots, \omega(e_{j_t}))\in {\bf B}\},$$
where ${\bf B}$ is an $i$-dimensional Borel set.  
If $k> m$, we simply define that
$${\cal F}_k={\cal F}_m.$$
Note that
$\cup_{i\leq k} W_i\subset \cup_{i\leq k+1} W_i$, so
$$ {\cal F}_0\subset {\cal F}_1\subset \cdots \subset {\cal F}_k \subset \cdots.$$
Thus, $\{{\cal F}_0, \cdots, {\cal F}_k,\cdots\}$ is a filtration constructed by random sets 
$\{W_1, \cdots, W_m\}$.
With these $\sigma$-fields, the martingale representation of
$T_n- E_\nu T_n$ is
$$T_n- E_\nu T_n=\sum_{k=1}^m \left[E_\nu(T_n\,\, |\,\, {\cal F}_k)- E_\nu (T_n\,\,|\,\, {\cal F}_{k-1})\right].\eqno{(3.2)}$$
Thus, 
$$ \{M_k= E_\nu(T_n\,\, |\,\, {\cal F}_k)-E_\nu T_n\}$$
defined an ${\cal F}_k$-martingale  sequence.
The martingale difference is defined by
$$\Delta_{k, n}= M_k-M_{k-1}.$$
It follows from the martingale property that
$$ T_n- E_\nu T_n=\sum_{k=1}^m \Delta_{k, n}\mbox{ and } 
\sigma_\nu (T_n)^2=\sum_{k=1}^m E_\nu \Delta_{k, m}^2.\eqno{(3.3)}$$

Let us carefully analyze the martingale difference $\{\Delta_{k, n}\}$.
By our definition, $E_\nu(T_n\,\, |\,\, {\cal F}_k)$ is ${\cal F}_k$ measurable.
For a set $\Gamma_k$, let
$\omega(\Gamma_k)$ and $c(\Gamma_k)$ be a random configuration set and a fixed configuration set, respectively.
 Thus,   we may regard 
$E_\nu(T_n\,\, |\,\, {\cal F}_k)$ as a random variable  on random set $\cup_{i\leq k} W_i=\Gamma_k$.
 For a random configuration  $\omega(\Gamma_k)=\omega(\cup_{i\leq k} W_i)$,
\begin{eqnarray*}
E_\nu(T_n\,\, |\,\, {\cal F}_k)&=&\sum_{c(B_e(n)\setminus \Gamma_k)}T_n(\omega(\Gamma_k), c(B_e(n)\setminus \Gamma_k))P_{\nu}(c(B_e(n)\setminus \Gamma_k))\\
&=&\sum_{c(B_e(n)\setminus \cup_{i\leq k} W_i)}T_n\left( 
 \omega(\bigcup_{i\leq k}W_i), c(B_e(n)\setminus  \bigcup_{i\leq k} W_i)\right)P_{\nu}\left( c(B_e(n)\setminus \bigcup_{i\leq k} W_i) \right),\hskip 1.2cm (3.4)
\end{eqnarray*}
where  $B_e$ is the edge set in $[-n, n]^2$ and 
$$P_{\nu}(c(B_e(n)\setminus \Gamma_k))=P_{\nu}(\omega(B_e(n)\setminus \Gamma_k)=c(B_e(n)\setminus \Gamma_k)),$$
and the sum takes all configurations $c(B_e(n)\setminus \bigcup_{i\leq k} W_i)$.
By this definition, we show the following lemma.\\

{\bf Lemma 3.1.} {\em For any   random set $\Gamma_k=\cup_{i\leq k} W_i$ constructed in  the spiral bubble structure, }
$$\sum_{c(B(n) \setminus \Gamma_{k})}P_\nu(  c(B(n)\setminus \Gamma_k))=1.\eqno{(3.5)}$$

{\bf Proof. } For  $\Gamma_k=\cup_{i\leq k} W_i$, let $I(\omega(\Gamma_k))$ be the indicator for 
the configuration of $\omega(\Gamma_k)$.  Thus,
$$E_\nu (E_\nu ( I(\omega(\Gamma_k))\,\,|\,\,{\cal F}_k))= E_{\nu} I(\omega(\Gamma_k)).\eqno{(3.6)}$$
Note that $I(\omega(\Gamma_k))$ is ${\cal F}_k$-measurable, so
if  we suppose that
$$\sum_{c(B(n) \setminus \Gamma_{k})}P_\nu( c(B(n)\setminus \Gamma_k))< 1,$$
then
$$E_\nu (E_\nu ( I(\omega(\Gamma_k))\,\,|\,\,{\cal F}_k))= E_\nu I(\omega(\Gamma_k)\left(
\sum_{c(B(n) \setminus \Gamma_{k})}P_\nu( c(B(n)\setminus \Gamma_k))\right)< E_{\nu} (I(\omega(\Gamma_k)).\eqno{(3.7)}
$$
The contradiction of (3.6) and (3.7) shows that Lemma 3.1 holds. $\Box$\\

\section{ Proof of  theorem.}
\subsection{A lower bound estimate of the theorem.}
For a given $\epsilon >0$, we may  decompose ${\cal H}_\epsilon$ defined in  Lemma 2.3 into ${\cal H}_\epsilon(k,t) \subset {\cal H}_\epsilon$ with
the first $k$, by following the order
${e_0, \cdots , e_m}$, such that
$D^*(v_k')$ contains its open cluster  $\bar{D}^*(v_k')$ with $t$ vertices larger than $ n^{2-5/4-\epsilon}$, where $v_k'$ and $v_k''$ are the two vertices of $e_k$. Thus,
$${\cal H}_\epsilon =\bigcup_{k} \bigcup_{t}{\cal H}_\epsilon(k,t).$$
By Lemma 2.3, note that $\{{\cal H}_\epsilon(k,t)\}$   are disjoint in $t$ and $k$, so there exists $C >0$ such that
$$P_\nu({\cal H}_\epsilon)=\sum_{k} \sum_{t}P_\nu ( {\cal H}_\epsilon(k,t)) > C.\eqno{(4.1)}$$

Now we choose a configuration with random set $\Gamma_k=\cup_{i\leq k} W_i$  in ${\cal H}_\epsilon(k,t)$
such that $e_k$ is open and $W_k=D^*(v_k')$, where $\{W_1, \cdots, W_m\}$ is constructed in the spiral bubble structure and $D^*(v_k')$ contains the  open cluster $\bar{D}^*(v_k')$ with vertex size $t$  larger than $ n^{2-5/48-\epsilon}$.
We denote by $I_t(\omega(\Gamma_k))$ the indicator of these configurations such that $e_k$ is open and the open cluster of $v'$ above has $t$ vertices   larger than
 $ n^{2-5/48-\epsilon}$. We also let $I(\Gamma_k,{\cal H}_\epsilon(k,t))$ be the indicator such that  on ${\cal H}_\epsilon (k, t)$ the first $k$ spiral bubble random sets
 are $\Gamma_k=\cup_{i\leq k} W_i$.
Thus,
$$I_t(\omega(\Gamma_k))\leq I(\Gamma_k, {\cal H}_\epsilon(k,t)).\eqno{(4.2)}$$
Indeed, if $I(\Gamma_k, {\cal H}_\epsilon(k,t))=1$, $e_k$ may or may not be open, but $e_k$ is open when
$I_t(\omega(\Gamma_k))=1$. In other words, if $I(\Gamma_k, {\cal H}_\epsilon(k,t))=1$, $W_k$ may be $W_k^+$ or $W_k^-$, but $W_k^+$ has to occur when $I_t(\omega(\Gamma_k))=1$. 

By (3.4), for  this random set  $\cup_{i\leq k-1} W_{i}\cup W_k^+$,
\begin{eqnarray*}
&&E_\nu(T_n\,\, |\,\, {\cal F}_{k})I(\Gamma_k, {\cal H}_\epsilon(k,t))\\
&=&\!\!\!\!\!\!\!\!\!\!\!\!\!\!\!\sum_{c(B_e(n)\setminus \cup_{i\leq k} W_i)}\!\!\!\!\!T_n\left( 
 \omega(\bigcup_{i\leq k-1}W_i), \omega(W_k^+), c(B_e(n)\setminus  \bigcup_{i\leq k} W_i)\right)P_{\nu}\left( c(B_e(n)\setminus \bigcup_{i\leq k} W_i) \right).\hskip 2.5cm (4.3)
 \end{eqnarray*}
In addition,  note that on $I(\Gamma_k, {\cal H}_\epsilon(k,t))$, $W_k$  can only take
$W_k^+$  or $W_k^-$, so by the second inequality in Lemma 2.3, 
\begin{eqnarray*}
&&
E_\nu(T_n\,\, |\,\, {\cal F}_{k-1})I(\Gamma_k, {\cal H}_\epsilon(k,t))\\
 &= & \!\!\!\!\!\sum_{c(B_e(n)\setminus \cup_{i\leq k} W_i)}\!\!\!T_n\left( 
 \omega(\bigcup_{i\leq k-1}W_i), \omega(W_k^+), c(B_e(n)\setminus  \bigcup_{i\leq k} W_i)\right)P_{\nu}\left( c(B_e(n)\setminus \bigcup_{i\leq k} W_i), c(W_k^+) \right)\\
 &+& \!\!\!\!\!\sum_{c(B_e(n)\setminus \cup_{i\leq k} W_i)}\!\!\!\!\! T_n\left( 
 \omega(\bigcup_{i\leq k-1}W_i), \omega(W_k^-), c(B_e(n)\setminus  \bigcup_{i\leq k} W_i)\right)P_{\nu}\left( c(B_e(n)\setminus \bigcup_{i\leq k}W_i) , c(W_k^-) \right)\\
 &\leq  & \!\!\!\!\!\sum_{c(B_e(n)\setminus \cup_{i\leq k} W_i)}\!\!\!T_n\left( 
 \omega(\bigcup_{i\leq k-1}W_i), \omega(W_k^+), c(B_e(n)\setminus  \bigcup_{i\leq k} W_i)\right)P_{\nu}\left( c(B_e(n)\setminus \bigcup_{i\leq k} W_i), \omega(e_k)=0 \right)\\
 &+& \!\!\!\!\!\sum_{c(B_e(n)\setminus \cup_{i\leq k} W_i)}\!\!\!\!\! T_n\left( 
 \omega(\bigcup_{i\leq k-1}W_i), \omega(W_k^-), c(B_e(n)\setminus  \bigcup_{i\leq k} W_i)\right)P_{\nu}\left( c(B_e(n)\setminus \bigcup_{i\leq k}W_i) , \omega(e_k)=1) \right)\\
 &=&0.5 \!\!\!\!\!\sum_{c(B_e(n)\setminus \cup_{i\leq k} W_i)}\!\!\!\!\!T_n\left( 
 \omega(\bigcup_{i\leq k-1}W_i), \omega(W_k^+), c(B_e(n)\setminus  \bigcup_{i\leq k} W_i)\right)P_{\nu}\left( c(B_e(n)\setminus \bigcup_{i\leq k} W_i) \right)\\
 &+&0.5 \!\!\!\!\!\!\!\sum_{c(B_e(n)\setminus \cup_{i\leq k} W_i)}\!\!\!\!\!T_n\left( 
 \omega(\bigcup_{i\leq k-1}W_i), \omega(W_k^-), c(B_e(n)\setminus  \bigcup_{i\leq k} W_i)\right)P_{\nu}\left( c(B_e(n)\setminus \bigcup_{i\leq k} W_i) \right).\hskip 2cm (4.4)
\end{eqnarray*}
If we put (4.3) and (4.4) together,  then for the random set $\Gamma_k=\bigcup_{i\leq k-1} W_i\cup W_k^+$ with $I_t(\omega(\Gamma_k))=1$, 
\begin{eqnarray*}
&&  \Delta_{k, n}\left(\omega(\bigcup_{i\leq k-1} W_i), \omega( W_k^+)\right)I(\Gamma_k, {\cal H}_\epsilon(k, t))\\
&\geq  &0.5 \!\!\!\!\!\!\!\!\!\sum_{c(B(n)\setminus \cup_{i\leq k} W_i)}\left [T_n\left(\omega(\bigcup_{i\leq k-1} W_i), \omega(W_k^+), c( B(n)\setminus \bigcup_{i\leq k} W_i)\right)
-T_n\left(\omega(\bigcup_{i\leq k-1} W_i), \omega(W_k^-), c( B(n)\setminus \bigcup_{i\leq k} W_i)\right)\right]\\
&&\times P_\nu\left( c( B(n)\setminus \bigcup_{i\leq k} W_i)\right). \hskip 11.5cm (4.5)
\end{eqnarray*}
Note that    $|\bar{D}^*(v_k')| =t\geq n^{2-5/28-\epsilon}$, so 
 $$
 T_n(\omega(\cup_{i\leq k-1} W_i), \omega(W_k^+), c( B(n)\setminus \cup_{i\leq k} W_i)
)-T_n(\omega(\cup_{i\leq k-1} W_i), \omega(W_k^-), c( B(n)\setminus \cup_{i\leq k} W_i))
\geq n^{2-5/48-\epsilon}.\eqno{(4.6)}$$
By (4.5) and (4.6), for the random set $\Gamma_k$ with $I_t(\omega(\Gamma_k))=1$, 
$$\Delta_{k, n}\left(\omega(\bigcup_{i\leq k-1} W_i), \omega( W_k^+)\right)I(\Gamma_k, {\cal H}_\epsilon(k, t))\geq 0.5 n^{2-5/48-\epsilon} \!\!\!\!\!\!\!\!\sum_{c(B(n)\setminus \cup_{i\leq k} W_i)}P_\nu\left(c( B(n)\setminus \bigcup_{i\leq k} W_i)\right).\eqno{(4.7)}$$
By  (4.7) and Lemma 3.1,  for each random set $\Gamma_k=\cup_{i\leq k} W_i\in {\cal H}_\epsilon(k,t)$ with $I_t(\omega(\Gamma_k))=1$,  there exists a constant $C >0$ such that
$$ \Delta_{k, n }^2 (\cup_{i\leq k-1} W_i\cup W_k^+)= \Delta_{k, n}^2(\omega(\Gamma_k))\geq C  n^{4-5/24-\epsilon}.\eqno{(4.8)}$$
Thus,   for each random set $\Gamma_k=\cup_{i\leq k} W_i\in {\cal H}_\epsilon(k,t)$ with $I_t(\omega(\Gamma_k))=1$,
 by (4.8), note that $\{I_t(\omega(\Gamma_k))\}$  are disjoint in $k, t, \Gamma_k,$
 and $\omega (\Gamma_k)$, so there exists $C>0$   such that
\begin{eqnarray*}
&&E_\nu \sum_{k} \Delta_{k, n}^2\geq E_\nu \sum_{k} \sum_{t}\sum_{\Gamma_k} \sum_{\omega(\Gamma_k)} I_t(\omega(\Gamma_k)) \Delta_{k, n}^2(\omega(\Gamma_k))\\
&\geq &  C n^{4-5/24-\epsilon} E_\nu \left(\sum_{k} \sum_{t}\sum_{\Gamma_k} \sum_{\omega(\Gamma_k)} I(\omega(\Gamma_k))\right). 
\hskip 8.5cm (4.9)
\end{eqnarray*}
By (4.9), 
$$E_\nu \sum_{k} \Delta_{k, n}^2\geq C n^{4-5/24-\epsilon}\sum_{k} \sum_{t} E_\nu I({\cal H}_\epsilon(k,t))I(\omega(e_k)=0).\eqno{(4.10)}$$
By (4.10) and Lemma 2.3,  there exist $C_i>0$  for $i=1,2$ such that 
$$E_\nu \sum_{k} \Delta_{k, n}^2\geq C_1 n^{4-5/24-\epsilon} P_\nu( {\cal H}_\epsilon) \geq 
C_2 n^{4-5/24-\epsilon}.\eqno{(4.11)}$$
By (3.3) and (4.11), note that $\epsilon$ can be arbitrarily small, so  the lower bound in the theorem  satisfies that
$$ \sigma^2( T_n)\geq n^{4-5/24+o(1)}.\eqno{(4.12)}$$
\subsection{ An upper bound estimate of the theorem.}
By (2.1), the upper bound in the theorem satisfies that
$$\sigma^2  (T_n)\leq E_\nu(T_n)^2 \leq C (E_\nu T_n)^2 \leq n^{4-5/24+o(1)}.\eqno{(4.13)}$$
\newpage
{\begin{center}{\large \bf References} \end{center}
Aizenman, M., Dulpantier, B. and Alharony, A. (1999). Path-crossing exponents and the external perimeter in 2D percolation. {\em Phy. Rev. Lett.}
{\bf 63}, 817--835.\\
Grimmett, G. (1999). {\em Percolation}. Springer-Verlag, New York.\\
Kesten, H. (1982). {\em Percolation Theory for Mathematicians.}  Birkhauser, Boston.\\
Kesten, H. (1986). The incipient infinite cluster in two-dimensional percolation
{\em Probab.  Theory  Related. Fields} {\bf 73}, 369--394.\\
Kesten, H. (1987). Scaling relations for 2D-percolation. {\em Comm.
Math. Phys.} {\bf 109}, 109--156.\\
Kesten, H., Sidoravicius,  V. and Zhang, Y. (1998). Almost all words are seen in critical site percolation on the triangular lattice. {\em Electron. J. Probab.} {\bf 3}  no. 10, 75pp.\\
Lawler, G. F., Schramm, O. and  Werner, W. (2002). One-arm exponent for critical 2D percolation. {\em Electron. J. Probab}. {\bf 7}, 13pp.\\
Smirnov, S. (2001). Critical percolation in the plane: conformal invariance, Cardy's formula, scaling limits. {\em  CR. Acad.  Sci. I-Math.} {\bf 333},  239--244.\\
Smirnov, S. and Werner, W. (2001). Critical exponent for two dimensional percolation. {\em Math.  Res. Lett. } {\bf 8}, 729--744.\\
Zhang, Y. (2001). A martingale approach in the study of percolation clusters on the ${\bf Z}^d$ lattice. {\em  J. Theoret Probab.} {\bf 14}, 165--187.\\

Yu Zhang\\
Department of Mathematics\\
University of Colorado, Colorado Springs\\
CO 80918 USA\\
yzhang3@uccs.edu

\end{document}